\documentclass{article}[15pt]
\usepackage{amsmath, epsfig,amssymb, multicol}
\usepackage[active]{srcltx}
\usepackage{epstopdf}
\usepackage{color}
\newtheorem{lemma}{Lemma}
\newtheorem{theorem}{Theorem}

\newtheorem{definition}{Definition}

\newtheorem{claim}{Claim}

\newcommand{\ex}{{\rm ex}}
\newcommand{\EX}{{\rm EX}}
\newcommand{\M}{{\cal M}}
\renewcommand{\l}{{\mathcal  L}}
\begin{document}
\title{Tur\'an number of the odd-ballooning of complete bipartite graphs}

\author{{Xing Peng and Mengjie Xia\footnote{Center for Pure Mathematics, School of Mathematical Sciences, Anhui University,
Hefei  230601, P.~R.~China. Email: {\tt x2peng@ahu.edu.cn}. Supported in part by the National Natural Science Foundation of China
(No.\ 12071002) and the Anhui Provincial Natural Science Foundation (No. 2208085J22). }}}

\date{}
\maketitle
\begin{abstract}
Given a graph $L$, the Tur\'an number  $\ex(n,L)$ is the maximum possible number of edges in an $n$-vertex $L$-free graph.
The study of Tur\'an number of graphs is a central topic in extremal graph theory.
Although the celebrated Erd\H{o}s-Stone-Simonovits theorem gives the asymptotic value of $\ex(n,L)$ for nonbipartite $L$, it is challenging in general to determine the exact value of $\ex(n,L)$ for $\chi(L) \geq 3$.
The odd-ballooning of $H$ is a graph such that each edge of $H$  is replaced by an odd cycle and all new vertices of  odd cycles are distinct. Here the length of odd cycles is not necessarily equal. The exact value of Tur\'an number of the odd-ballooning of $H$ is previously known for $H$ being a cycle, a path, a tree with assumptions, and $K_{2,3}$. In this paper, we manage to obtain the exact value of Tur\'an number of the odd-ballooning of $K_{s,t}$ with $2\leq s \leq t$, where  $(s,t) \not \in \{(2,2),(2,3)\} $ and each odd cycle has length at least five.




\end{abstract}

\section{Introduction}\label{sec-1}
Let $\l$ be a family of graphs. A graph $G$ is $\l$-free if $G$ does not contain any $L \in \l$ as a subgraph. The {\it Tur\'{a}n number} $\ex(n,\l)$ is the maximum number of edges in an $n$-vertex $\l$-free graph.     If an $n$-vertex $\l$-free graph $G$  contains exactly $\ex(n,\l)$ edges, then  $G$ is called an {\it extremal graph} for $\l$.  The set of extremal graphs with $n$ vertices is denoted by $\EX(n,\l)$.
For the case where $\l$ contains only one graph $L$, we write $\ex(n,L)$ and $\EX(n,L)$ for the Tur\'an number of $L$ and the set of extremal graphs for $L$ respectively. The well-known result by Mantel asserts that  $\ex(n,K_3)=\lfloor n^2/4\rfloor$ and the unique extremal graph is the balanced complete bipartite graph. Tur\'an's theorem \cite{Turan} gives the value of $\ex(n,K_{r+1})$ of all $ r \geq 2$. Moreover,
$\EX(n,K_{r+1})$ contains only the balanced complete $r$-partite graph which is known as {\it Tur\'an graph} $T_r(n)$.  The number of edges in the Tur\'an graph is denoted by $t_r(n)$.  Tur\'an's theorem is  considered as the origin of extremal graph theory. For the Tur\'an number of general graphs, Erd\H{o}s-Stone-Simonovits \cite{ES,ES1} theorem gives that
\[
\ex(n, L) = \left(1 - \frac{1}{\chi(L)-1} \right) \binom{n}{2}+o(n^2).
\]
Here $\chi(L)$ is the chromatic number of $L$.  This theorem is one of the cornerstones of extremal graph theory.

For $\chi(L) \geq 3$, the asymptotic value of $\ex(n, L)$ is given by Erd\H{o}s-Stone-Simonovits theorem. Thus it is meaningful to determine the exact value of $\ex(n, L)$ for such kind of graphs, which  is quite challenging in general. The notion of the decomposition family introduced by Simonovits \cite{S2} turns out to be very useful, for example \cite{chi,Liu,Yan,Yuan,Yuan1,Yuan2,Yuan3,ZC}. Given two graphs $G$ and $H$, the {\it join} $G \otimes H$ is the graph obtained from the vertex disjoint union of $G$ and $H$ by connecting each vertex of $G$ and each vertex of $H$.
\begin{definition}[Simonovits \cite{S2}]
Given a graph $L$ with $\chi(L)=p+1$,
the decomposition family $\M(L)$ is the set of minimal graphs $M$  such that $L\subset (M\cup \overline{K}_t)\otimes T_{p-1}((p-1)t)$, where    $t=t(L)$  is a constant.
\end{definition}
Roughly speaking,   $M \in \mathcal{M}(L)$ if the graph obtained from adding a copy of $M$ (but not any of its proper subgraphs) into a class of  $T_p(n)$ with $n$ large enough contains  $L$ as a subgraph. We remark  that the decomposition family $\mathcal{M}(L)$ always contains  bipartite graphs.
The following inequality can be checked easily:
\begin{equation}\label{lb1}
\ex(n,L) \geq t_p(n)+\ex(\lceil n/p\rceil, \M(L)).
\end{equation}
 To see this, let $G$ be a graph obtained from the Tur\'an graph $T_p(n)$ by adding a graph from $\EX(\lceil n/p\rceil, \M(L))$ to a largest part of $T_p(n)$. Thus $e(G)=t_p(n)+\ex(\lceil n/p\rceil, \M(L))$ and $G$ is $L$-free by the definition of $\M(L)$. Surprisingly, the above construction indeed gives the true value of $\ex(n,L)$ for many graphs, i.e.,
  $\ex(n,L)=t_p(n)+\ex(\lceil n/p\rceil, \M(L))$. Before we state examples, we introduce the following definitions and a related result.  A {\it matching} $M_k$ is a set of $k$ disjoint edges. For a graph $G$, the {\it matching number} $\nu(G)$ is the number of edges in a maximum matching of $G$. For positive integers $\nu$ and $\Delta$, let
  \[
  f(\nu,\Delta)=\max\left \{ e(G):\nu(G) \le \nu \hspace{0.5em} \mbox{and} \hspace{0.5em} \Delta(G) \le \Delta \right\}.
  \]
The following result   gives the upper bound for $f(\nu,\Delta)$.
\begin{theorem}[Chv\'{a}tal-Hanson \cite{CH}]\label{CH}
 $$f(\nu,\Delta)=\nu\Delta+\left \lfloor \frac{\Delta}{2}  \right \rfloor \left \lfloor \frac{\nu}{\left \lceil \Delta/2  \right \rceil }  \right \rfloor\le \nu\Delta+\nu.$$
 \end{theorem}
A special case where $\nu=\Delta=k-1$ was first proved by Abbott, Hanson, and Sauer \cite{AHS} as follows:
 \[
f(k-1, k-1)= \begin{cases}k^2-k & \text { if } k \text { is odd; } \\ k^2-\frac{3}{2} k & \text { if } k \text { is even. }\end{cases}
\]
 Let $F_k$ be the graph which consists of $k$ triangles sharing a common vertex.
 \begin{theorem}[Erd\H{o}s-F\"uredi-Gould-Gunderson  \cite{EFGG}] \label{EFGG}
 For $k \geq 1$ and $n \geq 50k^2$,
 \[
 \ex(n,F_k)=t_2(n)+f(k-1, k-1).
 \]
 \end{theorem}
One can see $\M(F_k)=\{M_k,K_{1,k}\}$, here  $K_{1,k}$ is a star with $k+1$ vertices. Moreover, $ \ex(\lceil n/2 \rceil, \{M_k,K_{1,k}\})=f(k-1, k-1)$. Thus the equality  holds in \eqref{lb1}
for $F_k$.

Given a graph $G$, let $kG$ denote the graph consisting of  $k$ vertex disjoint copies of $G$.
As $K_3$ is a clique,   Theorem \ref{EFGG} was generalized as follows.
\begin{theorem}[Chen-Gould-Pfender-Wei \cite{CGPW}]
For any $p \geq 2$ and $k \geq 1$, if $n \geq 16 k^3(p+1)^8$, then
$$
\operatorname{ex}\left(n, K_1 \otimes k K_p\right)=t_p(n)+f(k-1, k-1).
$$
\end{theorem}
It is clear that $F_k$ is a special case of $K_1 \otimes k K_p$ with $p=2$. One can observe that $\M(K_1 \otimes k K_p)=\{M_k,K_{1,k}\}$ and then   the equality  also holds in \eqref{lb1} for $K_1 \otimes k K_p$.
Motivated by this result, Liu \cite{Liu} introduced the concept of edge blow-up of graphs. For a graph $G$, the {\it edge blow-up} $G^{p+1}$ is a graph obtained from $G$ such that each edge is replaced by a clique with $p+1$ vertices and all new vertices for different cliques are distinct. One can see $F_k=K_{1,k}^{3}$ and $K_1 \otimes k K_p=K_{1,k}^{p+1}$.  So far, Tur\'an number of the edge blow-up of many families of graphs are known, for example, trees, cycles, keyrings, cliques $K_r$ with $p \geq r+1$, and complete bipartite graphs $K_{s,t}$ with $p \geq 3$,  see \cite{chi,Liu,NKSZ,WHLM,Yuan}. In general, a remarkable result by Yuan \cite{Yuan} gives the range of $\ex(n,G^{p+1})$ for $p \geq \chi(G)+1$.

If one view $K_3$ as  an odd cycle, then Theorem \ref{EFGG} can be generalized in another way. Let $G$ be a graph. The {\it odd ballooning} of $G$ is a graph obtained from $G$ such that each edge in $G$ is replaced by an odd cycle and all new vertices for different odd cycles are distinct. Apparently, $F_k$ is the odd-ballooning of $K_{1,k}$ in which all odd cycles are triangles. Notice that it is only meaningful to consider the odd-ballooning of bipartite graphs.  A star may be one of the simplest  bipartite graphs.
Hou, Qiu, and Liu \cite{HQL} first studied the Tur\'an number of the odd-ballooning of $K_{1,k}$ in which the length of all odd cycles is the same and is  at least five. Later, they \cite{HQL1} considered the general case in which triangles are allowed.  Zhu, Kang, and Shan \cite{ZKS} determined the Tur\'an number of the odd-ballooning of paths and cycles.
Recently, Zhu and Chen \cite{ZC} obtained the Tur\'an number of the odd-ballooning of  trees under some assumptions. It is great that previous results for paths and stars are special cases of  the result by Zhu and Chen \cite{ZC}.

As the Tur\'an number of the odd-ballooning of a cycle is known, the next step is to study  such a problem for bipartite graphs with many cycles. A possible candidate is complete bipartite graphs. Note that $K_{2,2}$ is $C_4$ and thus the simplest case is $K_{2,3}$. This case was solved by Yan \cite{Yan} previously.
The goal of this paper is to extend this result  to all complete bipartite graphs.

\begin{figure}[htbp]
\centerline{
\includegraphics[scale=0.8]{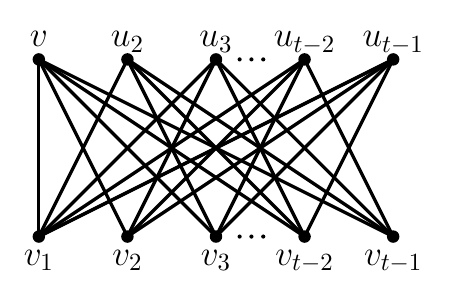}}
\caption{The graph $H$}
\label{fg1}
\end{figure}

To state our result, we need to define  a number of  graphs. Let $H$ be the graph obtained from $K_{t-1,t-1}$ by removing a $M_{t-2}$, see Figure \ref{fg1}, where $M_{t-2}=\{u_2v_2,\ldots,u_{t-1}v_{t-1}\}$. Note that $H$ is $P_4$ for $t=3$.
 For $2 \leq s \leq t$, let $G_{s,t}$ be the graph obtained from $T_2(n-s+1) \otimes K_{s-1}$ by embedding $H$ into one class of $T_2(n-s+1)$. Similarly, let $G_{3,3}'$ be the graph obtained from $T_2(n-2) \otimes K_{2}$  by embedding a triangle into one class of $T_2(n-2)$. The following theorem is our main result.
\begin{theorem}\label{main}
Let $F_{s,t}$ be the odd-ballooning of  $K_{s,t}$ with $t \geq s \geq 2$, where  $(s,t) \not \in \{(2,2),(2,3)\} $ and each  odd cycle has length at least 5. Then for $n$ large enough,

\[
\ex(n,F_{s,t})=\left \lceil \frac{n-s+1}{2} \right \rceil\left \lfloor \frac{n-s+1}{2} \right \rfloor+(s-1)(n-s+1)+\binom{s-1}{2}+t^{2}-3t+3.
\]
Moreover, $G_{s,t}$ is the only extremal graph for $t \geq 4$.  For $t=3$, there are at least two extremal graphs $G_{3,3}$ and $G'_{3,3}$.
\end{theorem}
The notation in this paper is standard. For a graph $G$ and a vertex $v \in V(G)$, let $N_G(v)=\{u: u \textrm{ is adjacent to } v\}$ be the {\it neighborhood} of $v$ and $d_G(v)=|N_G(v)|$ be the {\it degree} of $v$. If $X \subset V(G)$, then let  $d_G(v,X)=|N_G(v) \cap X|$ denote the number of neighbors of $v$ in $X$. Additionally, $G[X]$ is the subgraph induced by $X$, Let $e(X)$ and $\overline e(X)$ be the number of edges and non-edges  in $X$ respectively.  If $X$ and $Y$ are disjoint subsets of $V(G)$, then $e(X,Y)$ and $\overline e(X,Y)$ are the number of edges and non-edges between $X$ and $Y$ respectively.  If $e(X,Y)=|X||Y|$, then we say $X$ is completely adjacent to $Y$.
We use $uv$ to denote an  edge. If $u$ and $v$ are not adjacent, then we use $u \not \sim v$ to denote it.

The rest of this paper is organized as follows. In Section 2, we will recall a few results and prove some lemmas. In Section 3, we will present the proof of Theorem \ref{main}.

\section{Preliminaries}
The following  definition of symmetric graphs was introduced by Simonovits.

\begin{definition}[Simonovits \cite{S2}]
 Let $T_1$ and $T_2$ be connected subgraphs of $G$. They are called symmetric in $G$  if either $T_1=T_2$ or:

\noindent
 (1) $T_1 \cap T_2=\emptyset$; and

\noindent
(2) $(x, y) \notin G$ if $x \in T_1, y \in T_2$; and

\noindent
(3) there exists an isomorphism $\psi: T_1 \rightarrow T_2$ such that for every $x \in T_1$ and $u \in G-T_1-T_2, x$ is joined to $u$ if and only if $\psi(x)$ is joined to $u$.
\end{definition}
Note that $T_1, \ldots, T_k$ are symmetric if for every $1 \leq i<j \leq k,$ graphs  $T_i$ and $T_j$ are symmetric.
We also need to define a special family of  graphs.
\begin{definition}[Simonovits \cite{S2}]
 Let $\mathcal{D}(n, p, r)$ be the family of $n$-vertex graphs $G$ satisfying the following symmetric conditions:

 \noindent
(1) It is possible to omit at most $r$ vertices of $G$ such that the remaining graph $G^{\prime}$ is a join of graphs of almost equal order, i.e. $G^{\prime}=\otimes_{i=1}^p G^i$ where $\left|V\left(G^i\right)\right|=n_i$ and $\left|n_i-n / p\right| \leq r$ for any $i \in[p]$. The vertices in $V(G) \backslash V\left(G^{\prime}\right)$ are called the exceptional vertices.

\noindent
(2) For every $i \in[p]$, there exist connected graphs $H_i$ such that $G^i=k_i H_i$ where $k_i=n_i /\left|H_i\right|$ and $|V(H_i)| \leq r$. Moreover,  any two copies $H_i^j$ and $H_i^{\ell}$ in $G^i\left(1 \leq j<\ell \leq k_i\right)$ are symmetric subgraphs of $G$.
\end{definition}
 Our proof relies on the following theorem by  Simonovits.
\begin{theorem}[Simonovits \cite{S2}] \label{dnpr}
For a given graph $F$ with $\chi(F)=p+1 \geq 3$, if $\M(F)$ contains a linear forest, then there exist $r=r(F)$ and $n_0=n_0(r)$ such that $\mathcal{D}(n, p, r)$ contains an extremal graph for $F$ and $n \geq n_0$. Furthermore, if this is the only extremal graph in $\mathcal{D}(n, p, r)$, then it is the unique extremal graph for every sufficiently large $n$.
\end{theorem}

\begin{theorem}[Erd\H{o}s-Simonovits \cite{erdos67,erdos68,S1}]\label{min}
Let $F$ be a graph with chromatic number $t\geq 3$. If $G$ is an extremal graph for $F$ with $n$ large enough, then $\delta(G)= ((t-2)/(t-1) )n+o(n)$.
\end{theorem}


We conclude this section with two more lemmas.
\begin{lemma}\label{decomp}
Let $\mathcal{M}(F_{s,t})$ be the decomposition family of $F_{s,t}$. Then
$$
\left \{  K_{s,t},K_{1,p_1}\cup K_{1,p_2}\cup \cdots \cup K_{1,p_s}\cup M_q \right\} \subset  \mathcal{M}(F_{s,t}),
$$
where $\sum_{i=1}^s p_i+q=st$ and $0 \leq p_i \leq t$ for $1 \leq i \leq s$.
\end{lemma}
{\bf Proof:}  Let $A$ and $B$ be the two classes of $K_{s,t}$. For each $a \in A$ and $b \in B$, let $\ell_{a,b} \geq 5$ be the length of the odd cycle in $F_{s,t}$ associated with $a$ and $b$.  As a complete bipartite graph does not contain any odd cycles, if $M \in \mathcal{M}(F_{s,t})$, then $M$ contains at least one edge of each odd cycle and $|E(M)| \geq st$.
Given  $T_2(n)$ with  $n$ large, put $K_{1,p_1}\cup K_{1,p_2}\cup \cdots \cup K_{1,p_s}\cup M_q$ in one class.
View centers of $s$ stars as vertices of $A$ and put $B$ in the other class arbitrarily. Then for each $a \in A$ and $b \in B$, we can get an odd cycle $C_{a,b}$ either by including an edge incident with $a$ or by taking an edge from $M_q$. In addition, picking extra vertices appropriately, we can ensure that the length of $C_{a,b}$  equals $\ell_{a,b}$.  It is minimal as it contains $st$ edges. Similarly, we are able to show $K_{s,t} \in \M(F)$. \hfill $\square$

\vspace{0.2cm}
Note that if $p_i=0$ for each $1 \leq i \leq s$, then $K_{1,p_1}\cup K_{1,p_2}\cup \cdots \cup K_{1,p_s}\cup M_q$ is the matching $M_{st}$. Similarly, we can determine the decomposition family of $F_{1,t}$.
\begin{lemma}\label{decomp0}
Let $\mathcal{M}(F_{1,t})$ be the decomposition family of $F_{1,t}$. Then
\[
\mathcal{M}(F_{1,t})=\{K_{1,a} \cup M_{t-a}: 0 \leq a \leq t\}.
\]
\end{lemma}
{\bf Proof:} We have shown that each $M \in \mathcal{M}(F_{1,t})$ contains at least one edge of each odd cycle. Since $M$ is minimal, each odd cycle contributes exactly one edge to $E(M)$. Note that edges that are incident with the center of $K_{1,a}$ span a star and other edges form a matching. Meanwhile, one can find $F_{1,t}$ if one class  of $T_2(n)$ contains $K_{1,a} \cup M_{t-a}$. The lemma is proved. \hfill $\square$


\vspace{0.2cm}
\noindent
{\bf Remark 1:} From the proof of Lemma \ref{decomp}, one can see that the assumption $\ell_{a,b} \geq 5$ is crucial. Actually, if one can prove the theorem by assuming $\ell_{a,b}=5$ for each $a \in A$ and $b \in B$, then the proof also works for the case where $\ell_{a,b} \geq 5$  and the length of odd cycles is allowed to be different.

\section{Proof of Theorem \ref{main}}
\subsection{Proof of the lower bound}
We start to prove an auxiliary result. Recall $F_{s,t}$ is the odd-ballooning of  $K_{s,t}$ and  $G_{s,t}$ is the graph obtained from $T_2(n-s+1) \otimes K_{s-1}$ by embedding $H$ into one class of $T_2(n-s+1)$. In particular, $G_{1,t}$ is the graph obtained from $T_2(n)$ by embedding $H$ into one  class.
\begin{lemma} \label{base}
For any $1 \leq a  \leq (t+1)/2$, the graph $G_{1,t}$ does not contain $F_{a,t+1-a}$ as a subgraph, where all odd cycles have length at least five.
\end{lemma}
{\bf Proof:} Assume $V(K_{a,t+1-a})=A \cup B$, where $|A|=a$ and $|B|=t+1-a$.
Let $L \cup R$ be the partition of $G_{1,t}$ such that $H \subset L$. In addition, set $L'=\{v_1,\ldots,v_{t-1}\}$ and $R'=V(H) \setminus L'$. Suppose that $F_{a,t+1-a}$ is a subgraph of $G_{1,t}$ for some $a$.
Note that there are no odd cycles in a bipartite graph. Thus  each odd cycle in $F_{a,t+1-a}$   must contain an edge in $H$. In addition, all new vertices through the  operation of odd-ballooning are distinct. There are two cases.

\vspace{0.2cm}
\noindent
{\bf Case 1:} $A \cap L \neq \emptyset$.  If  $A \cap V(H) \neq \emptyset$, then we may assume $v_1 \in A$ for a moment.  As $A$ and $B$ are completely adjacent in $F_{a,t+1-a}$, one may assume $B=B' \cup B''$, where $B' \subset R'$ and $B'' \subset R$.

We begin with the case where both $B'$ and $B''$ are not empty. This implies  that $A \subset L'$, say $v_1,\ldots,v_a$.
Let $b'=|B'|$.
 As there are $t+1-a$ odd cycles associated with $v_1$ and each of them contains an edge in $H$, there is a  $(t+1-a-b')$-set $E_{v_1}$ of edges in $H$ that are associated with $v_1$.
 As $H$ is bipartite, let $B_1 \subset R' \setminus B'$ be the set which contains exactly one endpoint of edges in $E_{v_1}$. Note that vertices in $B_1$ are new in the operation of odd-ballooning.
  It is clear that $|B_1|=t+1-a-b'$. For each $2 \leq i \leq a$, if we repeat the analysis above for $v_i$, then there exists a subset $B_i \in R'\setminus B'$ associated with $v_i$ such that $B_i$ and $B_j$ are pairwise disjoint for $1 \leq i <j \leq a$. Here   vertices from $B_i$ and $B_j$  are new and of course are distinct.   Therefore,
  \begin{align*}
  t-1=|R'| &\geq |B' \cup B_1 \cup \cdots \cup B_a|\\
     &=|B'|+|B_1|+\cdots+|B_a|\\
     &=t+1-a+|B_2|+\cdots+|B_a| \\
     & \geq t+1-a+a-1=t.
  \end{align*}
  This is a contradiction.

  If  $B''$ is empty, then we consider locations  of vertices of $A$.  If $A \cap R=\emptyset$, then $A \subset L'$, say $A=\{v_1,\ldots,v_a\}$. Then $\{v_2,\ldots,v_a\}$ and $B$ form a $K_{a-1,t+1-a}$ which is a subgraph of $H-v_1$. The definition of $H$ gives that each $v_i$ (here $2 \leq i \leq t-1$) has a unique non-neighbor $u_i$ in $R'$. Furthermore, $u_i$ and $u_j$ are distinct. Therefore, vertices from $\{v_2,\ldots,v_a\}$ have at most $t-1-(a-1)=t-a$ common neighbors in $R'$, a contradiction.  If there is a vertex $r \in A \cap R$, here we note that $B \subset R'$, then there is a $(t+1-a)$-set $E_r$ of   edges from $H$ associated with  $r$ as there are $t+1-a$ odd cycles containing $r$ in $F_{a,t+1-a}$. Thus there is a $(t+1-a)$-subset $L_r' \subset L'$ which consists of exactly one endpoint of each edge from $E_r$. Note that vertices from $L_r'$ are new with respect to the operation of odd-ballooning. As the assumptions of $a$ and $|L'|=t-1$, there is only one such vertex. That is $A \setminus r \subset L'$. Note that $L_r'$ and $A \setminus r$ are disjoint. It holds that
  \[
  t-1=|L'| \geq |(A\setminus r) \cup L_r'|=|A\setminus r|+|L_r'|=a-1+t+1-a=t,
  \]
a contradiction.

If $B'$ is empty, then $A \subset L$. Reusing the argument above, for each vertex $a \in A$, there is a $(t+1-a)$-set $V_a \subset V(H)$ such that vertices from $V_a$ are new with respect to odd-ballooning. Thus $V_a$ and $V_{a'}$ are disjoint for $a \neq a' \in A$. Recall that we assume $v_1 \in A$. Thus $v_1 \cup_{a \in A} V_a \subset V(H)$.
It follows that $2(t-1) \geq a(t+1-a)+1$. Let $f(a)=-a^2+(t+1)a-2t+3$. Then $f(a)$ is concave down. As $a \leq (t+1)/2$, one can check $f(a)>0$ for $2 \leq a \leq (t+1)/2$, which is a contradiction to the inequality above. For the case of $a=1$, the graph $H$ must contain $M_p \cup K_{1,t-p}$ for some $0 \leq p \leq t$ as a subgraph by Lemma \ref{decomp0}, this is impossible by the definition of $H$.

To prove the case where  $A \cap V(H)\neq \emptyset$, it suffices to consider the case of $v_2 \in A$ by the symmetry of vertices in $H$. The proof  runs the same lines as above and it is skipped here.

To complete the proof for Case 1,  it remains to consider the case where  $A \cap V(H)=\emptyset$. Note that $A \subset L\setminus V(H) $ and  $B \subset R$. For each $a \in A$, if we let  $E_a \subset E(H)$ be the set of edges which consists of one edge from each odd cycle associated with $a$ in $F_{a,t+1-a}$, then $E_a$ is a matching with $t+1-a$ edges. Moreover, $E_a$ and $E_{a'}$ are disjoint for $a \neq a' \in A$. This is impossible by the definition of $H$. We complete the proof for Case 1.

\vspace{0.2cm}
\noindent
{\bf Case 2:} $A \cap L=\emptyset$. Equivalently, $A \subset R$  and $B \subset L$. By the symmetry of $A$ and $B$, this case can be proved by repeating the one for Case 1.  \hfill $\square$

\vspace{0.2cm}
 We are ready to prove the lower bound for $\ex(n,F)$.

\noindent
{\bf Proof of the lower bound:}   For positive integers $1 \leq k \leq t$,  we apply the induction on $k$ to show $G_{k,t}$ is $F_{k,t}$-free. The base case where $k=1$ follows from Lemma \ref{base} in which $a=1$.  Assume it is true for small $k$. To show the induction step, let $A$ and $B$ be the two classes of $K_{k,t}$, where $|A|=k$ and $|B|=t$. The vertex set of $K_{k-1}$ in $G_{k,t}$ is denoted by $W$. We view  $A$ and $B$ as subsets of $V(F_{k,t})$.   Suppose $F_{k,t}$ is a subgraph of $G_{k,t}$.
If $W \cap V(F_{k,t})=\emptyset$, then it means that $F_{k,t}$ is a subgraph of $G_{k,t}-W$. Note that $G_{k,t}-W$ is a subgraph of $G_{1,t}$ and $F_{1,t} \subset F_{k,t}$, then it must be the case where $F_{1,t} \subset G_{1,t}$. This is a contradiction by Lemma \ref{base}. Therefore, $W \cap V(F_{k,t}) \neq \emptyset$.

If there is a vertex $v \in V(F_{k,t})\setminus B$ such that $v \in W$, then  $F_{k,t}-v \subset G_{k,t}-v$.
Note that $F_{k,t}-v$ contains $F_{k-1,t}$ as a subgraph no matter $v \in A$ or $v$ is a new vertex with respect to odd-ballooning. Thus $G_{k,t}-v$ contains $F_{k-1,t}$ as a subgraph. However, this is impossible by the induction hypothesis and the observation $G_{k,t}-v \subset G_{k-1,t}$.

It is left to consider the case where $W \cap F_{k,t} \subset B$.
The assumption  $F_{k,t} \subset G_{k,t}$ implies that $F_{k,t}-W \subset G_{k,t}-W$.
 Note that $|W|=k-1$ and $W \cap F_{k,t} \subset B$.  It follows that $F_{k,t}-W$ contains $F_{k,t-k+1}$ as a subgraph.  Note that $G_{k,t}-W$ is a subgraph of $G_{1,t}$. We get that $F_{k,t-k+1}$ is a subgraph of $G_{1,t}$. This is a contradiction to  Lemma \ref{base}. There is a contradiction in each case and the proof for the induction step is complete. Therefore, $G_{s,t}$ does not contain $F_{s,t}$ as a subgraph.
As $G_n$ is an extremal graph, we get $e(G_n) \geq e(G_{s,t})$ and the lower bound for $e(G_n)$ follows.

%
%

%
%
%

\subsection{Proof of the upper bound}
Observe that $\chi(F_{s,t})=3$. Lemma \ref{decomp} gives that $\M(F_{s,t})$ contains a matching. By Theorem \ref{dnpr}, ${\cal D}(n,2,r)$ contains an extremal graph for $F_{s,t}$, say $G_n$,  provided $n$ is large enough. Here $r$ is a constant.

Assume $V(G_n)=A_1 \cup A_2 \cup R$, where $A_1$ and $A_2$ form a complete bipartite graph and $R$ is the set of exceptional vertices. By the definition of ${\cal D}(n,2,r)$, for $i \in \{1,2\}$, the subgraph $G[A_i]=k_iH_i$. If $H_i$ is nontrivial for some $i$, then $k_i \geq |A_i|/r \geq st$ as $|V(H_i)| \leq r$. This indicates that $M_{st} \subset G[A_i]$, a contradiction to Lemma \ref{decomp}. Therefore, both $A_1$ and $A_2$ are independent sets. Recall the definition of symmetric graphs. For each $1 \leq i \leq 2$ and each vertex $v \in R$, we get either $v$ is adjacent to all vertices in $A_i$ or $v$ has no neighbors in $A_i$. For $1 \leq i \leq 2$, we define
\[
B_i=\{v \in R: v \textrm{ is adjacent to all vertices in } A_{3-i}\}.
\]
Similarly, let
\[
W=\{v \in R: v \textrm{ is adjacent to all vertices in } A_1 \cup A_2\}
\]
and
\[
W'=\{v \in R: v \textrm{ has no neighbors in } A_1 \cup A_2\}.
\]
Apparently, $R=B_1 \cup B_2 \cup W \cup W'$ is a partition of $R$.  If $W'\neq \emptyset$, then for each vertex $v \in W'$, we have $d_{G_n}(v) \leq r$
 as $N_{G_n}(v) \subset R$ by the definition of $W'$. As $r$ is a constant, this is a contradiction to Theorem \ref{min}. Thus $W'=\emptyset$.

\begin{claim}
$|W|=s-1$.
\end{claim}
{\bf Proof:} Notice that  $W \cup A_1$ and $A_2$ form a complete bipartite graph. As $K_{s,t} \in \M(F_{s,t})$, we get $|W| \leq s-1$.  Suppose that $|W|=w \leq s-2$.
Note that
\[
e(G_n) \leq \left\lfloor \frac{n-w}{2} \right \rfloor \left\lceil \frac{n-w}{2} \right\rceil+e(B_1)+e(B_2)+\sum_{r \in W} d_{G_n}(r).
\]
Observe that
\[
\left\lfloor \frac{n-w}{2} \right \rfloor \left\lceil \frac{n-w}{2} \right \rceil
 \leq \left\lfloor \frac{n-s+1}{2} \right \rfloor + \left\lceil \frac{n-s+1}{2}\right\rceil+\frac{n(s-1-w)}{2},
\]
here adding one vertex to a balanced complete bipartite graph, one get at most $\tfrac{n}{2}$ new edges. Trivially,
\[
e(B_1)+e(B_2)+\sum_{r \in W} d_{G_n}(r) \leq \binom{|B_1|+|B_2|}{2}+(n-1)w \leq \binom{r}{2}+(n-1)w.
\]
Therefore, 
\begin{align*}
e(G_n) &\leq \left\lfloor \frac{n-s+1}{2} \right \rfloor + \left\lceil \frac{n-s+1}{2}\right\rceil+\frac{n(s-1-w)}{2}+\binom{r}{2}+(n-1)w \\
       &< \left\lfloor \frac{n-s+1}{2} \right \rfloor \left\lceil \frac{n-s+1}{2}\right\rceil +(s-1)(n-s+1)+\binom{s-1}{2}+t^2-3t+3,
\end{align*}
for the last inequality, we note that $r$ and $t$ are  constants, $w \leq s-2$, and $n$ is large enough.  This is a contradiction to the lower bound for $e(G_n)$ and the claim is proved. \hfill $\square$

\vspace{0.2cm}
Recall $F_{1,t}$ is the odd-ballooning of $K_{1,t}$. Let $T$ be the set of $t$ leaves of $K_{1,t}$.
\begin{claim}\label{base1}
The graph $G_n-W$ does not contain $F_{1,t}$ such that either $T \subset A_1$ or $T \subset A_2$.
\end{claim}
{\bf Proof:} Note that each vertex $w \in W$ is adjacent to all vertices of $A_1 \cup A_2$. Suppose that $G_n-W$ contains $F_{1,t}$ such that  $T \subset A_1$. Note that $W \cup T$ form a $K_{s-1,t}$.  As the definition of $W$, we get $F_{s-1,t}$ is a subgraph of $G_n$. Together with the $F_{1,t}$, there is an $F_{s,t}$ in $G_n$, a contradiction. \hfill $\square$

\vspace{0.2cm}
Let $\l=\{K_{1,p}\cup M_{t-p}: 0 \leq p \leq t\}$. To make the notation simple, we will use $B_1$ and $B_2$ to denote  subgraphs of $G_n$  induced by $B_1$ and $B_2$ respectively.
\begin{claim} \label{inside}
Both $B_1$ and $B_2$ are $\l$-free.
\end{claim}
{\bf Proof:} Suppose that $B_1$ contains $K_{1,p}\cup M_{t-p}$ for some $0 \leq p \leq t$. Consider the subgraph induced by $A_1\cup A_2 \cup B_1$. Observe that  $A_1 \cup B_1$ and $A_2$ form a complete bipartite graph. As $K_{1,p}\cup M_{t-p} \in \M(F_{1,t})$, the subgraph induced by $A_1\cup A_2 \cup B_1$ contains a copy of $F_{1,t}$, here we can choose $t$ vertices from $A_2$ as leaves of $K_{1,t}$, i.e., $T \subset A_2$. This is a contradiction by Claim \ref{base1}. We can show  $B_2$ is $\l$-free similarly. \hfill $\square$

\vspace{0.2cm}
Let $\nu_1=\nu(B_1)$ and $\nu_2=\nu(B_2)$.
\begin{claim}
$\nu_1+\nu_2 \leq t-1$.
\end{claim}
The proof is the same as the one for Claim \ref{inside} and it is skipped here.

\begin{claim}\label{nonneighbor}
 Assume that  $M_{p} \subset B_i$ and there is a vertex $v \in B_{3-i}$ with at least $q$ neighbors in $B_{3-i}$, where $i \in \{1,2\}$,  $1 \leq p,q \leq t-1$ and $p+q \geq t$. Then $v$ is not incident with at least $p+q-t+1$ edges from $M_p$.  Equivalently,   $v$ has at least $2(p+q-t+1)$ non-neighbors in $B_i$.
\end{claim}
{\bf Proof:} Suppose that $M_p \subset B_2$ and $v \in B_1$. Let $M_p=\{x_1y_1,\ldots,x_py_p\}$ and   $M'=\{x_iy_i: 1 \leq i \leq p, \textrm{   either } vx_i \in E(G_n) \textrm{ or } vy_i \in E(G_n)\}$. Note that $K_{1,q} \subset B_1$ with the center $v$.  We claim $|M'| \leq t-1-q$. Otherwise, suppose that $\{x_1y_1,\ldots,x_{t-q}y_{t-q}\} \subseteq M'$. Consider the subgraph induced by $A_1 \cup A_2 \cup B_1 \cup B_2$. Fix a $t$-subset $T$ of $A_2$. We can find the $F_{1,t}$ as follows. Actually, to get an odd cycle, we  include an edge from the star $K_{1,q}$ and  $M'$ one by one. For an edge $x_iy_i \in M'$, if  $vx_i$ is an edge and the length of the odd cycle associated with $x_iy_i$ is $2k+1$,  then we can find an odd cycle $vx_iy_ia_1a_2\cdots a_{2k-2}v$, here $a_{j} \in A_1$ for odd $j$,  $a_j \in A_2$ for even $j$, and $a_{2k-2} \in T$. This is a contradiction to Claim \ref{base1} and the claim is proved.  \hfill $\square$

\vspace{0.2cm}
Similarly, we can show the following variant and the proof is omitted here.
\begin{claim} \label{caset-1}
Let $v \in B_i$ and $\nu_i'$ be the matching number of the subgraph induced by  $B_i- (v \cup N_{B_i}(v))$. If $d_{B_i}(v)+\nu'_i=t-1$, then $v$ is not incident with any edge in $B_{3-i}$.
\end{claim}

\begin{claim} \label{lbound1}
$e(B_1 \cup B_2) \geq |B_1||B_2|+t^2-3t+3$.
\end{claim}
{\bf Proof:} Recall  the lower bound for $e(G_n)$. On the one hand,
$$
e(G_n-W) \geq \left\lfloor\frac{n-s+1}{2} \right\rfloor \left\lceil\frac{n-s+1}{2} \right\rceil+t^2-3t+3.
$$
On the other hand,
\begin{align*}
e(G_n-W)&=|A_1 \cup B_1||A_2|+|A_1||B_2|+e(B_1 \cup B_2)\\
        &=|A_1 \cup B_1||A_2 \cup B_2|-|B_1||B_2|+e(B_1 \cup B_2)\\
        &\leq \left\lfloor\frac{n-s+1}{2} \right\rfloor \left\lceil\frac{n-s+1}{2}\right \rceil-|B_1||B_2|+e(B_1 \cup B_2).
\end{align*}
The desired lower bound for $e(B_1 \cup B_2)$ follows.

\begin{claim}\label{degcase}
Assume that  both $B_1$ and $B_2$ contain edges. Let $\nu_1=\nu(B_1)$ and $\nu_2=\nu(B_2)$ with $\nu_1 \geq \nu_2$ and  $\nu_1+\nu_2 \leq t-1$. If either $(\nu_1,\nu_2)=(t-2,1)$ with $t \geq 5$ or $t \in \{3,4\}$, then $e(B_1\cup B_2) <|B_1||B_2|+t^2-3t+3$.
\end{claim}
{\bf Proof:}  Let $E'$ be the set of non-edges between $B_1$ and $B_2$.
The proof is split into the following cases.

\vspace{0.2cm}
\noindent
{\bf Case 1:} $(\nu_1,\nu_2)=(t-2,1)$ and $t \geq 5$.

If $\Delta(B_1)=t-1$, then let $v$ be  a vertex with the maximum degree and $N_{B_1}(v)=\{v_1,\ldots,v_{t-1}\}$. By Claim \ref{inside}, $B_1-\{v,v_1,\ldots,v_{t-1}\}$ is an independent set. Thus $e(B_1) \leq \sum_{i=1}^{t-1} d_{B_1}(v_i)$. If $d_{B_1}(v_i) \leq t-2$ for each $1 \leq i \leq t-1$, then $e(B_1) \leq t^2-3t+2$. As $\nu_2=1$,  edges in $B_2$ span a star or a triangle. Let $B_2' \subset B_2$  be the set of vertices with degree at least one.
 Then $|B_2'| \geq e(B_2)$ and the vertex $v$ is not adjacent to any vertex from $B_2'$ by Claim \ref{caset-1}, i.e., $|E'| \geq e(B_2)$. It is clear that $e(B_1\cup B_2) \leq e(B_1)+e(B_2)+|B_1||B_2|-|E'| < |B_1||B_2|+t^2-3t+3$. If there is a vertex $v_i$ such that  $d_{B_1}(v_i)=t-1$, then let $I=\{1 \leq i \leq t-1: d_{B_1}(v_i)=t-1\}$. As above, for each $i \in I$, the vertex $v_i$ is not adjacent to any vertex from $B_2'$. Thus
$|E'| \geq (|I|+1) e(B_2)$. Note that $e(B_1) \leq t^2-3t+2+|I|$ in this case. Therefore,
\begin{align*}
e(B_1\cup B_2) & \leq e(B_1)+e(B_2)+|B_1||B_2|-|E'| \\
                     &\leq |B_1||B_2|+t^2-3t+2+|I|+e(B_2)-(|I|+1)e(B_2)\\
                     &<|B_1||B_2|+t^2-3t+3.
\end{align*}
If $\Delta(B_1)=t-2$, then let $v$ be a vertex with the maximum degree and $N_{B_1}(v)=\{v_1,\ldots,v_{t-2}\}$, here $d_{B_1}(v_i) \leq t-2$ for each $1 \leq i \leq t-2$. Set $B_1'=B_1-(v \cup N_{B_1}(v))$. If $B_1'$ is an independent set, then
 $e(B_1) \leq \sum_{i=1}^{j} d_{B_1}(v_i) \leq t^2-4t+4$. If further $e(B_2) \leq t-2$, then the claim follows. Recall that $t \geq 5$ in this case. If $e(B_2)=t-1 \geq 4$, then edges in $B_2$ span a star since $\nu_2=1$. Let $u$ be the center of the star and then $u$ is not adjacent to vertices from $v \cup N_{B_1}(v)$, i.e., $|E'| \geq t-1$. Now, $e(B_1\cup B_2) \leq
|B_1||B_2|+t^2-4t+4+t-1-|E'|<|B_1||B_2|+t^2-3t+3$. If $B_1'$ contains edges, then $\nu(B_1')=1$ by Claim \ref{inside}.
 For $e(B_1') \leq t-2$, we get  $e(B_1) \leq \sum_{i=1}^{j} d_{B_1}(v_i)+e(B_1') \leq t^2-3t+2$. Note that there are at least $e(B_2)$ vertices which has degree at least one in $B_2$.  Claim \ref{caset-1} gives that $|E'| \geq e(B_2)$ and then $e(B_1 \cup B_2) \leq |B_1||B_2|+t^2-3t+2$. If $e(B_1')=t-1$, then $B_1'$ is a star as $t \geq 5$. This is impossible as it contradicts to the assumption $\Delta(B_1)=t-2$.

 If $\Delta(B_1) \leq t-3$, then $e(B_1) \leq t^2-4t+4$ by Theorem \ref{CH}. Repeating the argument above, we can show the upper bound for $e(B_1 \cup B_2)$.

\vspace{0.2cm}
\noindent
{\bf Case 2:} $t=4$. It is clear that either $(\nu_1,\nu_2)=(2,1)$ or $(\nu_1,\nu_2)=(1,1)$. For the case of $(\nu_1,\nu_2)=(1,1)$, note that $e(B_1) \leq 3$ and $e(B_2) \leq 3$. This implies that $e(B_1\cup B_2) \leq |B_1||B_2|+6<|B_1||B_2|+7$.  For the case of $(\nu_1,\nu_2)=(2,1)$, we can prove the upper bound for $e(B_1\cup B_2)$ by repeating the  argument  in Case 1.

\vspace{0.2cm}
\noindent
{\bf Case 3:} $t=3$. Apparently, $(\nu_1,\nu_2)=(1,1)$ in this case. Reusing the argument in Case 2, one can show the desired upper bound for $e(B_1\cup B_2)$. \hfill $\square$

\begin{claim}
One of $B_1$ and $B_2$ is an independent set.
\end{claim}
{\bf Proof:} Suppose that both $B_1$ and $B_2$ contain edges. Recall $\nu_1=\nu(B_1)$ and $\nu_2=\nu(B_2)$. We  assume that $\nu_1 \geq \nu_2 >0$. Claim \ref{inside} implies that $\Delta(B_1) \leq t-1$ and   $\Delta(B_2) \leq t-1$. If $\nu_1+\nu_2 \leq t-3$, then by Theorem \ref{CH}, we get
\[
e(B_1 \cup B_2) \leq |B_1|B_2|+ t \nu_1+t \nu_2 \leq |B_1||B_2|+ t^2-3t.
\]
A contradiction to Claim \ref{lbound1}. Thus it remains to consider the case where   $t-2 \leq \nu_1+\nu_2 \leq t-1$.  To simply the notation, let $K=G_n[B_1 \cup B_2]$ and show an upper bound for $e(K)$.  Observe that
 \[
 2e(K)=\sum_{v \in V(K)} d_{K}(v)=\sum_{v \in B_1}d_{K}(v)+ \sum_{v \in B_2}d_{K}(v).
 \]
We aim to establish $\sum_{v \in B_1}d_{K}(v) \leq |B_1||B_2|+2(t-\nu_2)\nu_1$. If $d_{K}(v,B_1)  \leq t-1-\nu_2$ for each vertex $v \in B_1$, then Theorem \ref{CH} gives that $e(B_1) \leq (t-\nu_2)\nu_1$, which yields that
\[\sum_{v \in B_1}d_{K}(v)=\sum_{v \in B_1} d_{K}(v,B_1)+ \sum_{v \in B_1} d_{K}(v,B_2) = 2e(B_1)+e(B_1,B_2) \leq 2(t-\nu_2)\nu_1+|B_1||B_2|.
\]
Thus we assume that there is a vertex $v \in B_1$ such that $d_{K}(v,B_1) \geq t-\nu_2$. Let $v_1$ be such a vertex with $d_{K}(v_1,B_1)=t-1+j_1-\nu_2$, here $j_1 \geq 1$. Applying Claim \ref{nonneighbor} with $i=1$, $p=\nu_2$, and $q=t-1+j_1-\nu_2$, we get that $v_1$ has at least $2j_1$ non-neighbors in $B_2$.
We change $j_1$ neighbors of $v_1$ from $B_1$ as non-neighbors arbitrarily and turn $2j_1$ non-neighbors of $v_1$ in $B_2$ as neighbors.
Let $K_1$ be the resulting graph. Observe that
$$
\sum_{v \in B_1} d_{K_1}(v)=\sum_{v \in B_1} d_{K}(v) \textrm{ and } d_{K_1}(v,B_2) \leq |B_2| \textrm{ for each } v \in B_1.
$$
Actually, changing $j_1$ neighbors of $v_1$ as non-neighbors makes the degree sum decreased by $2j_1$ while adding of $2j_1$ neighbors of $v_1$ contributes $2j_1$ to the degree sum.  Because of the choice of $v_1$, it is clear that $d_{K_1}(v,B_2) \leq |B_2|$ for each  $v \in B_1$. Moreover, $\nu(K_1[B_1]) \leq \nu(K[B_1])$.
For $i \geq 2$,
we shall define a vertex $v_i$ and a graph  $K_i$  recursively  such that $\sum_{v \in B_1} d_{K_i}(v)=\sum_{v \in B_1} d_{K_{i-1}}(v)$ and  $d_{K_i}(v,B_2) \leq |B_2|$ for each  $v \in B_1.$  If $d_{K_{i-1}}(v,B_1) \leq t-1-\nu_2$ for each $v \in B_1$, then we stop. Otherwise, let $v_i$ be a vertex such that $d_{K_{i-1}}(v_i,B_1)= t-1+j_i-\nu_2$ with $j_i \geq 1$. Observe that the added crossing edges so far are not incident with $v_i$. That is $N_{K_{i-1}}(v_i,B_2)=N_K(v_i,B_2)$.
Thus we can apply Claim \ref{nonneighbor} again to show that $v_i$ has at least $2j_i$ non-neighbors in $B_2$ (in $K_{i-1}$ and $K$). We repeat the operation above to get $K_i$ which satisfies desired properties.
 Assume that the process terminates after $\ell$ steps. We remark that all new crossing edges are distinct as they are associated with distinct vertices $v_1,\ldots,v_{\ell}$ in $B_1$.
Note that  $d_{K_{\ell}}(v,B_1)\leq  t-1-\nu_2$ for each $v \in B_1$ and $\nu(K_{\ell}[B_1]) \leq \nu(K[B_1])=\nu_1$. Therefore, by Theorem \ref{CH} and the definition of $K_i$ for $1 \leq i \leq \ell$, we get
\begin{align*}
\sum_{v \in B_1} d_{K}(v)&=\sum_{v \in B_1} d_K(v,B_1)+\sum_{v \in B_1} d_{K}(v,B_2)\\
                           &=\sum_{v \in B_1} d_{K_{\ell}}(v,B_1)+\sum_{v \in B_1} d_{K_{\ell}}(v,B_2)\\
                           & \leq \sum_{v \in B_1} d_{K_{\ell}}(v,B_1)+|B_1||B_2|\\
                           &=2e_{K_{\ell}}(B_1)+|B_1||B_2|\\
                           &\leq |B_1||B_2|+2(t-\nu_2)\nu_1.
\end{align*}
Repeating the argument above, we can show
$$
\sum_{v \in B_2} d_{K}(v) \leq |B_1||B_2|+2(t-\nu_1)\nu_2.
$$
Therefore,
\[
e(B_1\cup B_2)=E(K)=\frac{1}{2}\sum_{v \in B_1}d_{K}(v)+ \frac{1}{2}\sum_{v \in B_2}d_{K}(v)  \leq |B_1||B_2|+(t-\nu_2)\nu_1+(t-\nu_1)\nu_2.
\]
The case where either $t \in \{3,4\}$ or  $t\geq 5$ with $(\nu_1,\nu_2) = (t-2,1)$ is proved by Claim \ref{degcase}. We next assume $t \geq 5$ and $(\nu_1,\nu_2) \neq (t-2,1)$.
For the case of $\nu_1+\nu_2=t-1$, let $\nu_2=t-1-\nu_1$ and $g(\nu_1)=(t-\nu_2)\nu_1+(t-\nu_1)\nu_2=2\nu_1^2+(2-2t)\nu_1+t^2-t.$ Clearly, $g(\nu_1)$ is concave up. As assumptions
$\nu_1 \geq \nu_2$ and $(\nu_1,\nu_2) \neq (t-2,1)$,  the maximum value  of $g(\nu_1)$ is $g(t-3)=t^2-5t+12$.  Note that $g(t-3)<t^2-3t+3$ when $t \geq 5$.

For the case of $\nu_1+\nu_2=t-2$, let $\nu_2=t-2-\nu_1$ and $h(\nu_1)=2\nu_1^2+(4-2t)\nu_1+t^2-2t$. Similarly, the maximum value of $h(\nu_1)$ is $h(t-3)=t^2-4t+6<t^2-3t+3$ for $t \geq 4$.

We established $e(B_1 \cup B_2) < |B_1||B_2|+t^2-3t+3$ in each case. There is a contradiction to Claim \ref{lbound1} in each case and the proof is complete. \hfill $\square$

\vspace{0.2cm}
From now on,  assume $B_2$ is an independent set. Let $\Delta_1$ be the maximum degree of $B_1$. Recall $\nu_1$ is the matching number of $B_1$.
\begin{claim} \label{equal}
We have $e(B_1) \geq t^2-3t+3$.
Furthermore, if $e(B_1) \leq t^2-3t+3+k$ for an integer $0 \leq k \leq s-1$, then $\overline e(W)+\overline e(W,B_1) \leq k$. As a special case where $k=0$, i.e., $e(B_1)=t^2-3t+3$,  then $e(B_1,B_2)=|B_1||B_2|$ and each vertex $w \in W$ has degree $n-1$.
\end{claim}
{\bf Proof:} As $G_n$ is an extremal graph, it satisfies that
$$
e(G_n) \geq  \left\lfloor\frac{n-s+1}{2} \right\rfloor \left\lceil\frac{n-s+1}{2} \right\rceil+(s-1)(n-s+1)+\binom{s-1}{2}+t^2-3t+3.
$$
Let $C_1=A_1 \cup B_1$ and $C_2=A_2 \cup B_2$. Observe that
\[
e(G_n)=e(C_1, C_2)+e(B_1)+e(W,C_1 \cup C_2)+e(W).
\]
Notice that $e(C_1, C_2) \leq \left\lfloor\tfrac{n-s+1}{2} \right\rfloor \left\lceil\tfrac{n-s+1}{2} \right\rceil$ and
$e(W,C_1 \cup C_2)+e(W) \leq (s-1)(n-s+1)+\binom{s-1}{2}$.
Therefore, $e(B_1) \geq t^3-3t+3$.
If $e(B_1) \leq t^2-3t+3+k$ for an integer $0 \leq k \leq s-1$, then
\[
e(G_n) \leq \left\lfloor\frac{n-s+1}{2} \right\rfloor \left\lceil\frac{n-s+1}{2} \right\rceil+t^2-3t+3+k+(s-1)(n-s+1)+\binom{s-1}{2}-\overline e(W)-\overline e(W,B_1).
\]
Combined with the lower bound for $e(G_n)$, we get $\overline e(W)+\overline e(W,B_1)\leq k$ and the second part follows. The special case of $k=0$ can be shown similarly.  \hfill $\square$

\begin{claim} \label{mdeg}
The subgraph $B_1$ has exactly one connected component and $\Delta_1=t-1$.
\end{claim}
{\bf Proof:}  Claim \ref{inside} gives that $\Delta_1 \leq t-1$ and $\nu_1 \leq t-1$.  We assert that $\Delta_1 \geq t-2$. Otherwise, $e(B_1) \leq t^2-3t+2$ by Theorem \ref{CH}. This is a contradiction to Claim \ref{equal} and the assertion follows. Let $v$ be a vertex with maximum degree in $B_1$ and $C_1$ be the connected component containing $v$.

If $\Delta_1=t-1$, then  $C_1$ is the only connected component in   $B_1$. Otherwise, $K_{1,t-1} \cup M_1 \subset B_1$. This is a contradiction to Claim \ref{inside}.
If $\Delta_1=t-2$, then let $N(v)$ be the neighborhood of $v$ in $B_1$.
Suppose $B_1$ has another connected component $C_2$. Then $\nu(C_2)=1$. Otherwise, $K_{1,t-2} \cup M_2 \subset B_1$, which is a contradiction to Claim \ref{inside}. Similarly, we can show $B_1-(v \cup N(v))$ is an independent set, i.e., each edge in $C_1$ is incident with a vertex from $N(v)$. As $\Delta_1=t-2$, we get  $e(C_2) \leq t-2$ for $t \geq 5$. Therefore, $e(B_1) \leq \sum_{i=1}^{t-2} d_{B_1}(v_i)+e(C_2) \leq (t-2)^2+t-2=t^2-3t+2$ provided $t \geq 5$.
For $t=4$, if $e(C_2) \leq 2$, then $e(B_1) \leq 6$ which leads to the same contradiction. If $t=4$ and $e(C_2)=3$, then $C_2$ is a triangle. If $v_1v_2$ is  an edge, then $e(C_1)=3$ and $e(B_1)=6$, a contradiction to Claim \ref{equal}. If $v_1$ is not adjacent to $v_2$, then it must be the case where $e(C_1)=4$ and $e(C_2)=3$ as the lower bound for $e(B_1)$.
 Thus we get $\nu(C_1)=2$ and $\Delta(C_2)=2$ which imply that $K_{1,2} \cup M_2 \subset B_1$. This is a contradiction to Claim \ref{inside}.
For $t=3$, as $\Delta_1=1$, we get $e(B_1)=\nu(B_1)=2<3$ and obtain the same contradiction.  Therefore, $B_1$ has exactly one connected component.

Notice that we already showed $\Delta_1 \geq t-2$ and $B_1$ contains exactly one connected component. One can repeat the argument above to show $\Delta_1$ must equal $t-1$. \hfill $\square$

\begin{claim}\label{st3}
If $t=3$, then the unique connected component in $B_1$  either is a triangle or is a $P_4$.
\end{claim}
{\bf Proof:} Note that $s=t=3$ now. By Claim \ref{mdeg}, there is a vertex $v$ with two neighbors in $B_1$, say $v_1$ and $v_2$.   Claim \ref{inside} yields that $B_1 \setminus (N_{B_1}(v_1) \cup N_{B_1}(v_2))$ is an independent set.
If $v_1v_2$ is an edge, then $vv_1v_2$ is a triangle and there is no other edge as Claim \ref{mdeg}. If $v_1 \not \sim v_2$, then both $v_1$ and $v_2$ have at most one  neighbor other than $v$ as $\Delta_1=2$. If one of $v_1$ and $v_2$ has degree one, then we get a $P_4$ and we are done. Assume $v_1u_1$ and $v_2u_2$ are two edges. If $u_1\neq u_2$, then $K_{1,2} \cup M_1 \subset B_1$, this is a contradiction to Claim \ref{inside}. If $u_1=u_2$, then $e(B_1)=4$. Note that $s-1=2$ and assume $W=\{w_1,w_2\}$.
If both $w_1$ and $w_2$ are completely adjacent to $B_1$, then $\{v,u_1,w_1\}$ and $\{v_1,v_2,w_2\}$   form a $K_{3,3}$. Otherwise, as the lower bound and the upper bound for $e(G_n)$,  only one of $w_1$ and $w_2$ has a unique non-neighbor in $B_1$, say $w_1$.  Thus $w_1$ is completely adjacent to one of $\{v_1,v_2\}$ and $\{v,u_1\}$, say $\{v,u_1\}$. Note that $w_2$ is  completely adjacent to $B_1$.
 Now $\{v,u_1,w_2\}$ and $\{v_1,v_2,w_1\}$   form a $K_{3,3}$. Thus $F_{3,3} \subset  G_n$ by Lemma \ref{decomp}. This is a contradiction and the claim is proved. \hfill $\square$

\begin{claim} \label{nonadj}
Let $u$ and $v$ be two vertices from $B_1$ such that $u$ and $v$ are not adjacent. If $d_{B_1}(u)=t-1$, then $N_{B_1}(v) \subseteq N_{B_1}(u)$.
\end{claim}
{\bf Proof:}  Suppose that $v$ has a neighbor $x$ such that $x \not \in N_{B_1}(u)$. Then  observe that $K_{1,t-1} \cup M_1 \subset B_1$, here $K_{1,t-1}$ is formed by $u$ and its neighbors in $B_1$ while $M_1$ is the edge $vx$. This is a contradiction to Claim \ref{inside}. \hfill $\square$

\begin{claim} \label{indep}
Let $v$ be a vertex from $B_1$ with $d_{B_1}(v)=t-1$. Then $N_{B_1}(v)$ is an independent set for $t \geq 4$.
\end{claim}
 {\bf Proof:}  Let $N_{B_1}(v)=\{v_1,\ldots,v_{t-1}\}$ and  $B_1'=B_1-(v \cup N_{B_1}(v))$. Without causing any confusion, for each $u \in B_1$,  we will use $N(u)$ to denote $N_{B_1}(u)$ in the proof.  Claim \ref{inside} gives us that $B_1'$ is an independent set. Thus
\begin{equation}\label{count}
e(B_1)=\sum_{i=1}^{t-1} d_{B_1}(v_i)-e(N(v)).
\end{equation}
This equation will be used frequently in the proof.
Suppose $N(v)$ contains an edge.   Then there is a vertex $v_i \in N(v)$ such that $d_{B_1}(v_i)=t-1$. Otherwise, the equation \eqref{count} gives that $e(B_1) \leq (t-1)(t-2)-1=t^2-3t+2$ and this is a contradiction to Claim \ref{equal}. Without loss of generality, let $v_1$ be such a vertex.

 The vertex $v_1$ is adjacent to some $v_j \in N(v)$. If it is not the case, then $N(v_1)=v \cup T$, where $T \subset B_1'$ and $|T|=t-2$. By Claim \ref{nonadj}, $N(v_i) \subseteq N(v_1)=v \cup T$ for each $2 \leq i \leq t-1$  and then $N(v)$ is an independent set, a contradiction to the assumption. Let $\{v_2,\ldots,v_{p}\}$ be  the set of vertices that are adjacent to $v_1$. We next show  $p<t-1$. Otherwise,  $e(N(v)) \geq t-2$ as $p=t-1$. By equation \eqref{count} and Claim \ref{equal},
\[
t^2-3t+3 \leq e(B_1) \leq (t-1)^2-e(N(v)) \leq t^2-3t+3,
\]
which implies that $d_{B_1}(v_i)=t-1$ for each $1 \leq i \leq t-1$, $e(B_1)=t^2-3t+3$, and $e(N(v))=t-2$, i.e., $\{v_2,\ldots,v_{t-1}\}$ is an independent set. By Claim \ref{nonadj}, there is a  subset $D \subset B_1'$ with $t-3$ vertices  such that $N(v_i)=\{v,v_1\} \cup D$ for each $2 \leq i \leq t-1$. Note that $N(v_1)=\{v,v_2,\ldots,v_{t-1}\}$.
If we let $L=\{v_2,\ldots,v_{t-1}\}$ and $R=\{v,v_1\} \cup D$, then $L \cup R$ form a $K_{t-2,t-1}$.
As $e(B_1)=t^2-3t+3$, each vertex from $W$ has degree $n-1$ by Claim \ref{equal}. Thus including vertices from $W$ properly, there is a $K_{s,t}$ in $B_1\cup W$. As $A_1 \cup B_1 \cup W$ and $A_2$ form a complete bipartite graph, then $F_{s,t} \subset G_n$ by Lemma \ref{decomp}, a contradiction. Thus $p<t-1$.


There is a vertex $v_i$ such that $p+1 \leq i \leq t-1$ and $d_{B_1}(v_i)=t-1$. If there is no such a vertex, then by equation \eqref{count}
\begin{align*}
e(B_1)&=\sum_{i=1}^{t-1} d_{B_1}(v_i)-e(N(v))\\
       &=\sum_{i=1}^{p} d_{B_1}(v_i)+\sum_{i=p+1}^{t-1} d_{B_1}(v_i)-e(N(v))\\
       & \leq p(t-1)+(t-1-p)(t-2)-e(N(v))\\
       & \leq p(t-1)+(t-1-p)(t-2)-(p-1)\\
       &=t^2-3t+3,
\end{align*}
here $e(N(v)) \geq p-1$ as $v_1$ is adjacent to $v_i$ for each $1 \leq i \leq p$ . Together with Claim \ref{equal}, we get

\noindent
(1) $d_{B_1}(v_i)=t-1$ for $1 \leq i \leq p$,

\noindent
(2) $d_{B_1}(v_i)=t-2$ for each $p+1 \leq i \leq t-1$,

\noindent
(3) $e(N(v))=p-1$, i.e., $E(N(v))=\{v_1v_2,\ldots,v_1v_p\}$.

Assume $N(v_1)=\{v,v_2,\ldots,v_p\} \cup T$, where $T \subset B_1'$ with $|T|=t-p-1$.
Claim \ref{nonadj} implies that $N(v_i) \subset v \cup T$ for each $p+1 \leq i \leq t-1$.
Since $d_{B_1}(v_i)=t-2$ and $|v \cup T| = t-p \leq t-2$, we get $p=2$ and $N(v_i) =v \cup T$ for each $3 \leq i \leq t-1$.
Similarly, observe that  $N(v_2)=\{v,v_1\} \cup T$. Now $v \cup T$ and $\{v_1,v_2,\ldots,v_{t-1}\}$ form a $K_{t-2,t-1}$.
Notice that each vertex from $W$ has degree $n-1$ as $e(B_1)=t^2-3t+3$. We can show $F_{s,t} \subset G_n$ similarly and this is a contradiction.
In the following, we assume that $d_{B_1}(v_i)=t-1$ for $p+1 \leq i \leq p+q$.

We next show $p=2$.   Recall that  $N(v_1)=\{v,v_2,\ldots,v_p\} \cup T$, where $T \subset B_1'$ with $|T|=t-p-1$.
As $v_1 \not \sim v_i$ for $p+1 \leq i \leq p+q$, Claim \ref{nonadj} gives that $N(v_i) \subseteq N(v_1)=\{v,v_2,\ldots,v_p\}  \cup T$. The assumption  $d_{B_1}(v_1)=d_{B_1}(v_i)=t-1$ indicates that $v_i$ is adjacent to all vertices from $\{v_2,\ldots,v_p\}$ for each $p+1 \leq i \leq p+q$.  Thus $e(N(v)) \geq p-1+q(p-1)$. If $p>2$, then by equation \eqref{count}, we get $e(B_1) \leq (t-1)(p+q)+(t-1-p-q)(t-2)-(p-1)(q+1)=t^2-3t+3+(2-p)q<t^2-3t+3$, which is a contradiction to Claim \ref{equal}. Therefore, $p=2$ and $e(B_1)=t^2-3t+3$. This implies that $d_G(v_i)=t-2$ for $q+3 \leq i \leq t-1$ and each $w \in W$ has degree $n-1$ by Claim \ref{equal}. Moreover, $E(N(v))=\{v_1v_2,v_3v_2,\ldots,v_{q+2}v_2\}$.
 Notice that $|T|=t-3$ and  $N(v_1)=\{v,v_2\} \cup T$.  As $v_i \not \sim v_1$ and $d_{G_n}(v_i)=t-1$ for each $3 \leq i \leq q+2$,   we get $T \subset N(v_i)$ by Claim \ref{nonadj}. Note that $\{v_1,v_3,\ldots,v_{q+2}\}$ and $\{v,v_2\} \cup T$ form a $K_{q+1,t-1}$.  As $q+t+s-1 \geq s+t$ and each $w \in W$ has degree $n-1$,   then    $\{v_1,v_3,\ldots,v_{q+2}\} \cup \{v,v_2\} \cup T \cup W$
 contains a $K_{s,t}$ by including vertices from $W$ appropriately.
 Since $A_1 \cup B_1 \cup W$ and $A_2$ is a complete bipartite graph, $F_{s,t} \subset G_n$ by Lemma \ref{decomp}. This is a contradiction and the claim is proved.
\hfill $\square$


\begin{claim} \label{oneH}
The subgraph $B_1$ contains exactly one copy of $H$ as a subgraph for $t \geq 4$.
\end{claim}
{\bf Proof:} We reuse the assumptions in the proof of Claim \ref{indep}. Then $v$ is a vertex from $B_1$ with maximum degree, $N_{B_1}(v)=\{v_1,\ldots,v_{t-1}\}$, and  $B_1'=B_1-(v \cup N_{B_1}(v))$. Claim \ref{indep} tells us that $N_{B_1}(v)$ is an independent set.
Notice that $e(B_1)=\sum_{i=1}^{t-1} d_{B_1}(v_i)$. Claim \ref{equal} implies that there is at least one $v_i \in N(v)$ such $d_{B_1}(v_i)=t-1$. Assume that $d_{B_1}(v_i)=t-1$ for each $1 \leq i \leq k$. We claim $k=1$. Otherwise,
Claim \ref{nonadj} yields that $N(v_i)=v \cup T$ for each $1 \leq i \leq k$, where $T \subset B_1'$ and $|T|=t-2$. Then $L=\{v_1,\ldots,v_k\}$ and $R=v \cup T$ form a $K_{k,t-1}$ with $k \geq 2$. In the following, we will show $K_{s,t} \subset B_1 \cup W$. As $A_1 \cup B_1 \cup W$ and $A_2$ is a complete bipartite graph, then $F_{s,t}$ is a subgraph of  $G_n$ by Lemma \ref{decomp} which is a contradiction to the assumption.
There are two cases depending on $k$.

\vspace{0.2cm}
\noindent
{\bf Case 1:}  $2 \leq k \leq s-1$.  Let
$$
W'=\{w \in W: w \textrm{ is adjacent to all vertices in }  R\}.
$$
Note that each $w \in W \setminus W'$ has a non-neighbor in $B_1$. As $e(B_1) \leq t^2-3t+2+k$, it follows that $|W\setminus W'| \leq  \overline e(W,B_1)+\overline e(W) \leq k-1$ by Claim \ref{equal}, i.e., $|W'| \geq s-k$. Assume $|W'|=s-k+j$.

For $j=0$,  each vertex from $W\setminus W'$ has exactly one non-neighbor in $R$ and has degree  $n-2$.
Note that $|W'| \leq s-2$. Pick $w \in W\setminus W'$ arbitrarily and notice that $w$ is adjacent to all vertices in $L \cup W'$. Now, $L \cup W'$ and $R \cup w$ form a $K_{s,t}$.

For $j>0$, we have $|W\setminus W'|=k-j-1$.
 As $\overline e(W \setminus W',B_1) \geq k-j-1$, we get $\overline e(W')+\overline e(W',L) \leq j$ by Claim \ref{equal}. If $\overline e(W',L)=0$, then there is a vertex $w \in W'$ which has at least $s-k$ neighbors in $W'$ provided $s-k+j>2$. Otherwise, $e(W') \leq \tfrac{(s-k-1)(s-k+j)}{2}<\binom{s-k+j}{2}-j$, a contradiction. Let $w \in W'$ be such a vertex and $W''$ be the set of $s-k$ neighbors of $w$ in $W'$.  Then $L \cup W''$ and $R \cup w$ is a $K_{s,t}$.
If $s-k=j=1$, then $k=s-1$, $|W'|=2$, and $|W\setminus W'|=s-3$. Assume $W'=\{w_1,w_2\}$. If $w_1$ is adjacent to $w_2$, then $L \cup w_1$ and $R \cup w_2$ is a $K_{s,t}$. Thus we assume $w_1$ is not adjacent to $w_2$.  For $W\setminus W' \neq \emptyset$, then each $w \in W\setminus W'$ is completely adjacent to $L \cup W'$.
  We can find a $K_{s,t}$ in $B_1 \cup W$ as we did in the case where $j=0$. For $W\setminus W' = \emptyset$, i.e., $s=3$ and $W=W'=\{w_1,w_2\}$. As $w_1$ is not adjacent to $w_2$, we get both $w_1$ and $w_2$ are completely adjacent to $L \cup R$ by Claim \ref{equal}.
  Note that $t>3$ by the assumption. Then $d_{G_n}(v_i)=t-2$ for $3 \leq i \leq t-1$ and $v_i$ has a unique non-neighbor in $T$  by Claim \ref{nonadj}. Let $Y=\{v_3,\ldots,v_{t-2}\}$.
    Thus there are two vertices $t_1,t_2 \in T$ such that both $t_1$ and $t_2$ are completely adjacent to $Y$ as $|T|=t-2$. Now $\{v,t_1,t_2\}$ and $W \cup \{v_1,v_2\} \cup Y$ form a $K_{3,t}$.

If $\overline e(W',L)>0$, then let $Z \subset W'$ be the set of vertices which have non-neighbors in $L$. Thus $\overline e(W') \leq j-|Z|$. Recall $|W'|=s-k+j$.
Removing vertices from $W'$ which has non-neighbors in $W'$  one by one, we will get a subset $W''$ such that $|W''| \geq s-k+|Z|$ and  $W''$ is a clique. Let $w \in W'' \setminus Z$ and $W''' \subset W''$ be the set of $s-k$ neighbors of $w$ in $W''$. Then $L \cup W'''$ and $R \cup w$ form a $K_{s,t}$.

\vspace{0.2cm}
\noindent
{\bf Case 2:}  $k \geq s$. We first consider the case where $k \geq s+1$. We claim that there is a vertex $w \in W$ such that $w$ has at least $s$ neighbors in $L$. If there is a such vertex $w$, then let $L'$ be the set of $s$ neighbors of $w$ in $L$. Observe that $L'$ and $R \cup w$ is $K_{s,t}$.  It is left to show the existence of $w$.  Let
$$W'=\{w \in W: w \textrm{ has  at most  } s-1 \textrm{  neighbors in }  L\}.$$
  As $e(B_1) \leq t^2-3t+2+k$ and each $w \in W'$ has at least $k-s+1$ non-neighbors in $L$,  it follows that $(k-s+1)|W'| \leq k-1$. Thus $|W'| \leq 1+\tfrac{s-2}{k-s+1} \leq 1+\tfrac{s-2}{2}$ since $k \geq s+1$.
  If $s \geq 4$, then $|W'| \leq s-2$ and the desired vertex $w \in W$ exists. If  $s=3$, then the inequality above gives $|W'| \leq 1$ and the existence of $w$ also follows. For $s=2$, note that $W$ contains only one vertex, say $w$,  and $e(B_1) \leq t^2-3t+2+k$. Notice that $w$ has at most one neighbor in $L$, i.e.,  $d_{G_n}(w) \leq n-k$,
    Combining the upper bound and the lower bound for $e(G_n)$, we get that $d_{G_n}(w) = n-k$ and $d_{G_n}(v_i)=t-2$ for each $k+1 \leq i \leq t-1$, here we assume $k<t-1$ for a while.  Claim \ref{nonadj} give that $N_{G_n}(v_i)=v \cup T_i$ such that $T_i \subset T$ and $|T_i|=t-3$, i.e., $v_i$ has exactly one non-neighbor in $T$. As $|T|=t-2$ and $k \geq s+1$, there is a vertex $u_1 \in T$ such that $u_1$ is adjacent to each vertex $v_i$ for $k+1 \leq i \leq t-1$.  As $d_{G_n}(w) = n-k$  and $w$ has at least $k-1$ non-neighbors in $L$, then $w$ is completely adjacent to $R$. Now $\{v,u_1\}$ and $w \cup N_{B_1}(v)$ form a $K_{2,t}$.
For $k=t-1$,  note that $R$ and $L \cup w$ is a $K_{t-1,t}$ and we can find a $K_{s,t}$ easily.

It remains to prove the case of $k=s$. As above, if $|W'| \leq s-2$, then we are done. Thus we assume $|W'|=s-1$, i.e., $W'=W$. Similarly, one can show $d_{G_n}(w)=n-2$ for each $w \in W$ and the unique non-neighbor of $w$ is in $L$.  In addition, $d_{G_n}(v_i)=t-2$ for each $s+1 \leq i \leq t-1$, here we also assume $s<t-1$ for a moment. As above,  we can find an $(s-1)$-subset $T' \subset T$ such that  $T'$ and $\{v_{s+1},\ldots,v_{t-1}\}$ form a complete bipartite graph. Now $v\cup T'$ and $w \cup N_{B_1}(v)$ is a $K_{s,t}$ for any $w \in W$.
For $k=s=t-1$, as above, $R$ and $L \cup w$ is a $K_{t-1,t}$ for any $w \in W$ and we are able to find $K_{s,t}$ easily.

There is a contradiction in each case and then $k=1$ follows.  Assume $v_1$ is the unique vertex such that $d_{B_1}(v_1)=t-1$ and $N(v_1)=v \cup T$, where $T \subset B_1'$ with $|T|=t-2$. Meanwhile,  $e(B_1)=t^2-3t+3$ and $d_{B_1}(v_i)=t-2$ for each $2 \leq i \leq t-1$. By Claim \ref{nonadj}, $N_{B_1}(v_i) \subset v \cup T$ for each $2 \leq i \leq t-1$. Actually, for each $v_i$, there is a unique vertex $u_i \in T$ such that $v_i \not \sim u_i$. Moreover, we assert that $u_i \neq u_j$ for $2 \leq i \neq j \leq t-2$. Otherwise, suppose $u_2=u_3$ without losing any generality. Observe that  $\{v_1,v_2,v_3\}$ and $v \cup (T \setminus u_2)$ is  $K_{3,t-2}$. Claim \ref{equal} together with $e(B_1)=t^2-3t+3$ imply that each vertex from $W$ has degree $n-1$. Let $W'$ be an $(s-3)$-subset of $W$.
  Then $\{v_1,v_2,v_3\} \cup W'$ and $v \cup (T \setminus u_2) \cup (W\setminus W')$ form a $K_{s,t}$ for $s \geq 3$. For $s=2$, there is a vertex $u_1 \in T$ which is completely adjacent to $\{v_2,\ldots,v_{t-1}\}$. Now $\{v,u_1\}$ and $N(v) \cup W$ is a $K_{2,t}$.
By Lemma \ref{decomp}, $F_{2,t} \subset G_n$. This is a contradiction and the proof is complete.

\vspace{0.2cm}
\noindent
{\bf Proof of Theorem \ref{main}:} For $(s,t) \neq (3,3)$, Claim \ref{oneH} gives that $B_1$ contains exactly one copy of $H$ as a subgraph. Let $C_1=A_1 \cup B_1$ and $C_2=A_2 \cup B_2$.  Reusing the proof for Claim \ref{lbound1}, we get that
\[
e(G_n) \leq \left\lfloor\frac{n-s+1}{2} \right\rfloor \left\lceil\frac{n-s+1}{2} \right\rceil+(n-s+1)(s-1)+\binom{s-1}{2}+t^2-3t+3.
\]
As the lower bound for $e(G_n)$, it follows that $d_{G_n}(w)=n-1$ for each $w \in W$. In addition, $e(C_1,C_2)=\left\lfloor\frac{n-s+1}{2} \right\rfloor \left\lceil\frac{n-s+1}{2} \right\rceil$, i.e., $C_1$ and $C_2$ is a balanced complete bipartite graph with $n-s+1$ vertices. Therefore, $G_n=G_{s,t}$ and it is the unique extremal graph  by Theorem \ref{dnpr}. For $(s,t)=(3,3)$, one can show either $G_n=G_{3,3}$ or $G_n=G_{3,3}'$. For this case, we are not able to characterize all extremal graphs. \hfill $\square$


\begin{thebibliography}{99}
\bibitem{AHS}
H.~Abbott, D.~Hanson, and H.~Sauer, Intersection theorems for systems of sets, {\it J.~Combin.~Theory Ser.~A,} {\bf 12} (1972), 381--389.
\bibitem{CGPW}
G.~Chen, R.~Gould, F.~Pfender, and B.~Wei, Extremal graphs for intersecting cliques, {\it J.~Combin.~Theory Ser.~B,} {\bf 89} (2003), 159--171.
 \bibitem{chi}
  C.~Chi and L.~Yuan, The Tur\'{a}n number for the edge blow-up of trees: The missing case, {\it Discrete Math.,} {\bf 346}(6) (2023), No.113370.
\bibitem{CH}
V.~Chv\'{a}tal and D.~Hanson, Degrees and matchings, {\it J.~Combin.~Theory Ser.~B,} {\bf 20} (1976), 128-138.
\bibitem{erdos67}
 P.~Erd\H{o}s, Some recent results on extremal problem in graph theory, Theory of Graphs (ed P. Rosenstiehl), (Internat. Sympos., Rome, 1966), Gordon and Breach, New York, and Dunod, Paris, 1967, 117--123.
\bibitem{erdos68}
P.~Erd\H{o}s, On some new inequalities concering extremal properties of graphs, Theory of Graphs (P. Erd\H{o}s and G. Katona, Eds.), Academic Press, New. York, 1968, 77--81.
\bibitem{ES}
P.~Erd\H{o}s and M.~Simonovits, A limit theorem in graph theory, {\it Studia Sci.~Math Hungar.,} {\bf 1} (1966), 51--57.
\bibitem{ES1} P.~Erd\H{o}s and A.~Stone, On the structure of linear graphs, {\it Bull.~Amer.~Math.,} {\bf 52} (1946), 1089--1091.
\bibitem{EFGG}
P.~Erd\H{o}s, Z.~F\"{u}redi, R.~Gould, and D.~Gunderson, Extremal graphs for intersecting triangles, {\it J.~Combin.~Theory Ser.~B,} {\bf 64}(1) (1995), 89--100.

\bibitem{HQL}
X.~Hou, Y.~Qiu, and B.~Liu, Extremal graph for intersecting odd cycles, {\it Electron.~J. Combin.,} {\bf 23}(2) (2016), P29.
\bibitem{HQL1}
X.~Hou, Y.~Qiu, and B.~Liu, Tur\'an number and decomposition number of intersecting odd cycles, {\it Discrete Math.,} {\bf 341}(1) (2018), 126--137.

\bibitem{Liu}
 H.~Liu, Extremal graphs for blow-ups of cycles and trees, {\it Electron.~J.~Combin.,} {\bf 20}(1) (2013), P65.
\bibitem{NKSZ}
Z.~Ni, L~Kang, E.~Shan, and H.~Zhu, Extremal graphs for blow-ups of keyrings, {\it Graphs Combin.,} {\bf 36}(6) (2020),  1827--1853.
\bibitem{S1}
M.~Simonovits,  A method for solving extremal problems in graph theory, stability problems, in:
Theory of Graphs, Proc.~Colloq., Tihany, 1966, Academic Press, New York, 1968, 279--319.
\bibitem{S2}
M.~Simonovits, Extremal graph problems with symmetrical extremal graphs, additional chromatic conditions, {\it Discrete Math.,} {\bf 7} (1974), 349--376.
\bibitem{Turan}
P.~Tur\'{a}n, On an extremal problem in graph theory (in Hungrarian), {\it Mat.~es Fiz.~Lapok.,} {\bf 48} (1941), 436--452.
\bibitem{WHLM}
A.~Wang, X~Hou, B.~Liu, and Y.~Ma,
The Tur\'an number for the edge blow-up of trees,
{\it Discrete Math.,} {\bf 344}(12) (2021), No.112627.
\bibitem{Yan}
N.~Yan, The Tur\'an number of graphs with given decomposition family, {\it Acta Scientiarum Naturalium Universitatis Nankaiensis,} {\bf 54}(4) 2021, 34--43.
\bibitem{Yuan3} L.~Yuan, Extremal graphs for the $k$-flower, {\it J.~Graph Theory,} {\bf 89}(1) (2018), 26--39.
\bibitem{Yuan2} L.~Yuan, Extremal graphs for odd wheels, {\it J.~Graph Theory,} {\bf 98}(4) (2021),  691--707.
\bibitem{Yuan} L.~Yuan, Extremal graphs for edge blow-up of graphs, {\it J.~Combin.~Theory Ser.~B,} {\bf 152}  (2022), 379--398.
\bibitem{Yuan1}
 L.~Yuan, Extremal graphs of the $p$th power of paths, {\it European J.~Combin.,} {\bf 104} (2022),  No.103548.

\bibitem{ZKS}
H.~Zhu, L.~Kang, and E.~Shan,  Extremal Graphs for odd-ballooning of paths and cycles, {\it Graphs Combin.,} {\bf 36}(3) (2020),  755--766.
\bibitem{ZC}
X.~Zhu and Y.~Chen, {\rm Tur\'{a}n} number for odd-ballooning of trees, {\it J.~Graph Theory,} {\bf 104}(2) (2023), 261--274.
\end{thebibliography}
\end{document}